# ON ECCENTRIC CONNECTIVITY INDEX

**Bo Zhou** and **Zhibin Du**

Department of Mathematics, South China Normal University,

Guangzhou 510631, China

e-mail: zhoubo@scnu.edu.cn; zhibindu@126.com

**Abstract**

The eccentric connectivity index, proposed by Sharma, Goswami and Madan, has been employed successfully for the development of numerous mathematical models for the prediction of biological activities of diverse nature. We now report mathematical properties of the eccentric connectivity index. We establish various lower and upper bounds for the eccentric connectivity index in terms of other graph invariants including the number of vertices, the number of edges, the degree distance and the first Zagreb index. We determine the $n$-vertex trees of diameter $d$, $3 \leq d \leq n-2$, with the minimum eccentric connectivity index, and the $n$-vertex trees of $p$ pendent vertices, $3 \leq p \leq n-2$, with the maximum eccentric connectivity index. We also determine the $n$-vertex trees with respectively the minimum, second-minimum and third-minimum, and the maximum, second-maximum and third-maximum eccentric connectivity indices for $n \geq 6$.

## 1. INTRODUCTION

A topological index is a numerical descriptor of the molecular structure derived from the corresponding (hydrogen-depleted) molecular graph. Various topological indices are widely used for quantitative structure-property relationship (QSPR) and quantitative structure-activity relationship (QSAR) studies [1-4].

Let $G$ be a simple connected graph with vertex set $V(G)$. For a vertex $u \in V(G)$, $d_G(u)$ or $d_u$ denotes the degree of $u$, $e_G(u)$ or $e_u$ denotes the eccentricity of $u$ in $G$. Recall that $d_u$ is the number of (first) neighbors of $u$ in $G$, while $e_u$ is the path length from $u$ to a vertex $v$ that is farthest from $u$, i.e., $e_u = \max\{d(u,v\,|\,G) : v \in V(G)\}$, where $d(u,v\,|\,G)$ denotes the distance between $u$ and $v$ in $G$.

Sharma, Goswami and Madan [5] introduced an adjacency-cum-distance based topological index, the eccentric connectivity index, of the graph $G$, defined as

$$\xi^c = \xi^c(G) = \sum_{u \in V(G)} d_u e_u,$$

which has been employed successfully for the development of numerous mathematical models for the prediction of biological activities of diverse nature [5-13], and modified forms have also been given, for example in [14] by taking into consideration the presence as well as relative position of the heteroatom(s) in a molecular structure. So it is of interest to study the mathematical properties of this invariant.

In this paper, we give some basic mathematical properties of the eccentric connectivity index. We establish various lower and upper bounds in terms of other graph invariants including the number of vertices, the number of edges, the degree distance and the first Zagreb index. We determine the $n$-vertex trees of diameter $d$, $3 \le d \le n-2$, with the minimum eccentric connectivity index, and the $n$-vertex trees with respectively the minimum, second-minimum and third-minimum eccentric connectivity indices for $n \ge 6$. We also determine the $n$-vertex trees of $p$ pendent vertices, $3 \le p \le n-2$, with the maximum eccentric connectivity index, and the $n$-vertex trees with respectively the maximum, second-maximum and third-maximum eccentric connectivity indices for $n \ge 6$. It appears that the eccentric connectivity index satisfies the basic requirement to be a branching index.

## 2. PRELIMINARIES

For a connected graph $G$, the radius $r(G)$ and diameter $D(G)$ are, respectively, the minimum and maximum eccentricity among vertices of $G$. A connected graph is called a self-centered graph if all of its vertices have the same eccentricity [15]. Then a connected graph $G$ is self-centered if and only if $r(G) = D(G)$.

Let $P_n$ and $S_n$ be respectively the path and the star with $n$ vertices. Let $K_n$ be the complete graph with $n$ vertices. Let $C_n$ be the cycle with $n \geq 3$ vertices. Let $K_{r,s}$ be the complete bipartite graph with $r$ vertices in one vertex-class and $s$ vertices in the other vertex-class.

By direct calculation, the following formulae hold: $\xi^c(K_n) = n(n-1)$, $\xi^c(K_{r,s}) = 4rs$ for $r, s \geq 2$, $\xi^c(S_n) = 3(n-1)$ for $n \geq 3$, $\xi^c(C_n) = 2n\left\lfloor \dfrac{n}{2} \right\rfloor$, $\xi^c(P_n) = \left\lfloor \dfrac{3(n-1)^2+1}{2} \right\rfloor$. For $P_n$, we have two cases based on the parity of $n$. If $n \geq 2$ is even, then

$$\xi^c(P_n) = 2\left(n-1+\sum_{i=n/2}^{n-2} 2i\right) = 4\sum_{i=n/2}^{n-1} i - 2(n-1)$$
$$= n\left(\dfrac{n}{2}+n-1\right) - 2(n-1) = \dfrac{3(n-1)^2+1}{2},$$

and if $n \geq 3$ is odd, then

$$\xi^c(P_n) = 2\left(n-1+\sum_{i=(n+1)/2}^{n-2} 2i\right) + 2 \cdot \dfrac{n-1}{2} = 4\sum_{i=(n+1)/2}^{n-1} i - (n-1)$$
$$= (n-1)\left(\dfrac{n+1}{2}+n-1\right) - (n-1) = \dfrac{3(n-1)^2}{2}.$$

## 3. RESULTS FOR GENERAL GRAPHS

In this section, we give lower and upper bounds for the eccentric connectivity index of connected graphs in terms of graph invariants such as the number of vertices, the number of edges, the radius, the diameter, the degree distance and the first Zagreb index.

**Proposition 1.** *Let $G$ be a connected graph with $m$ edges. Then*

$$2mr(G) \leq \xi^c(G) \leq 2mD(G)$$

*with either equality if and only if $G$ is a self-centered graph.*

**Proof.** Note that $\sum_{u \in V(G)} d_u = 2m$. It is easily seen that $\xi^c(G) = \sum_{u \in V(G)} d_u e_u \leq \sum_{u \in V(G)} d_u D(G)$
$= 2mD(G)$ with equality if and only if $e_u = D(G)$ for any $u \in V(G)$, i.e., $G$ is a self-centered graph. The lower bound follows similarly by using $e_u \geq r(G)$ for $u \in V(G)$. □

Let $G$ be a connected graph with $m$ edges. If $r(G) \geq 2$, then by Proposition 1, $\xi^c(G) \geq 4m$ with equality if and only if $G$ is a self-centered graph of radius two. The self-centered graphs with radius two of minimum size have been characterized [15, Theorem 2.7].

**Corollary 1.** *Let $G$ be a connected graph with $n \geq 4$ vertices for which the complement $\overline{G}$ is also connected. Then*

$$\xi^c(G) + \xi^c(\overline{G}) \geq 2n(n-1)$$

*with equality if and only if both $G$ and $\overline{G}$ are self-centered graphs with radius two.*

**Proof.** Let $m$ and $\overline{m}$ be respectively the number of edges of $G$ and $\overline{G}$. Evidently, $2(m+\overline{m}) = n(n-1)$. Since both $G$ and $\overline{G}$ are connected, each has radius at least two, and then by Proposition 1,

$$\xi^c(G) + \xi^c(\overline{G}) \geq 4m + 4\overline{m} = 2n(n-1)$$

with equality if and only if $G$ and $\overline{G}$ are self-centered, and $r(G) = r(\overline{G}) = 2$. □

If $G$ is a connected graph with $n$ vertices and $m$ edges, and the number of vertices of degree $n-1$ in $G$ is equal to $k$, where $0 \leq k \leq n$, then it is easily seen that $\xi^c(G) \geq (n-1)k + 2[2m - k(n-1)] = 4m - k(n-1)$ with equality if and only if all the $n-k$ vertices of degree less than $n-1$ have eccentricity two. We use this fact in the proofs of Propositions 2 and 3.

A 3-vertex connected graph is either $S_3$ or $K_3$, for which $\xi^c(S_3) = \xi^c(K_3) = 3(n-1)$.

**Proposition 2.** *Let $G$ be an $n$-vertex connected graph, where $n \geq 4$. Then $\xi^c(G) \geq 3(n-1)$ with equality if and only if $G = S_n$.*

**Proof.** Let $m$ be the number of edges, and $k$ the number of vertices of degree $n-1$ in $G$, where $0 \leq k \leq n$. Then $\xi^c(G) \geq 4m - k(n-1)$ with equality if and only if all the $n-k$ vertices of degree less than $n-1$ have eccentricity two.

If $k = 0$, then $\xi^c(G) \geq 2mr(G) \geq 4m \geq 4(n-1) > 3(n-1)$. Suppose that $k \geq 1$. Note that all vertices of $G$ except the ones with degree $n-1$ are of degree at least $k$. We have $2m \geq k(n-1) + k(n-k)$, and thus,

$$\xi^c(G) \geq 2[k(n-1) + k(n-k)] - k(n-1) = k(3n - 2k - 1).$$

Obviously, the function $f(x) = x(3n - 2x - 1)$ with $1 \leq x \leq n$ attains the minimum value for $x = 1$ or $n$. Note that $f(1) = 3(n-1) < n(n-1) = f(n)$. Hence, $\xi^c(G) \geq f(1) = 3(n-1)$ with equality if and only if $k = 1$ and $m = n-1$, i.e., $G = S_n$. □

Denote by $G \vee H$ the graph formed from vertex-disjoint graphs $G$ and $H$ by adding edges between each vertex in $G$ and each vertex in $H$. For positive integers $n$ and $m$ with $n - 1 \leq m < \binom{n}{2}$, let $a = a_{n,m} = \left\lfloor \dfrac{2n - 1 - \sqrt{(2n-1)^2 - 8m}}{2} \right\rfloor$. Then $a < n$. Let $\mathbf{G}_{(n,m)}$ be the set of graphs $K_a \vee H$, where $H$ is a graph with $n - a$ vertices and $m - \binom{a}{2} - a(n-a)$ edges. Note that $a$ is the largest integer satisfying $2m \geq a(n-1) + a(n-a)$, i.e., $h(a) \geq 0$ with $h(a) = a^2 - 2na + a + 2m$, we have $\left[ m - \binom{a}{2} - a(n-a) \right] - (n - a - 1) = \dfrac{1}{2} h(a+1) < 0$. Thus, each vertex of $H$ has eccentricity two in $K_a \vee H$.

**Proposition 3.** *Let $G$ be an n-vertex connected graph with $m$ edges, where $n-1 \le m < \binom{n}{2}$. Let $a = \left\lfloor \dfrac{2n-1-\sqrt{(2n-1)^2 - 8m}}{2} \right\rfloor$. Then*

$$\xi^c(G) \ge 4m - a(n-1)$$

*with equality if and only if $G \in \mathbf{G}_{(n,m)}$.*

**Proof.** It may be easily checked that $a \ge 1$. Let $k$ be the number of vertices of degree $n-1$ in $G$, where $0 \le k \le n-1$. If $k = 0$, then $\xi^c(G) \ge 4m > 4m - a(n-1)$. Suppose that $k \ge 1$. Then $\xi^c(G) \ge 4m - k(n-1)$ with equality if and only if all the $n-k$ vertices of degree less than $n-1$ have eccentricity two. Note that $2m \ge k(n-1) + k(n-k)$ and $a$ is the largest integer satisfying $2m \ge a(n-1) + a(n-a)$. We have $k \le a$, and thus $\xi^c(G) \ge 4m - a(n-1)$ with equality if and only if $G$ has exactly $a$ vertices of degree $n-1$ and all other vertices have eccentricity two, i.e., $G \in \mathbf{G}_{(n,m)}$. □

Note that for $n \ge 4$, $a_{n,n} = 1$ and $\mathbf{G}_{(n,n)}$ contains exactly the unicyclic graph formed by adding an edge to the star $S_n$, and for $n \ge 5$, $a_{n,n+1} = 1$ and $\mathbf{G}_{(n,n+1)}$ contains exactly two bicyclic graphs formed by adding two edges to the star $S_n$. Thus, by Proposition 3, we have

**Corollary 2.** *Let $G$ be a unicyclic graph with $n \ge 4$ vertices. Then*

$$\xi^c(G) \ge 3n + 1$$

*with equality if and only if $G$ is formed by adding an edge to the star $S_n$.*

**Corollary 3.** *Let $G$ be a bicyclic graph with $n \ge 5$ vertices. Then*

$$\xi^c(G) \ge 3n + 5$$

with equality if and only if $G$ is formed by adding two edges to the star $S_n$.

Corollaries 2 and 3 also follow from Proposition 1: For a unicyclic graph $G$ with $n \geq 4$ vertices, if $r(G) \geq 2$, then $\xi^c(G) \geq 4m = 4n > 3n+1$, and if $r(G) = 1$, then $G$ is formed by adding an edge to the star $S_n$, for which $\xi^c(G) = (n-1) \cdot 1 + 2 \cdot 2 + 2 \cdot 2 + 1 \cdot 2 \cdot (n-3) = 3n+1$. This proves Corollary 2. For a bicyclic graph $G$ with $n \geq 5$ vertices, if $r(G) \geq 2$, then $\xi^c(G) \geq 4m = 4(n+1) > 3n+5$, and if $r(G) = 1$, then $G$ is formed by adding two edges to the star $S_n$, for which $\xi^c(G) = (n-1) \cdot 1 + 2 \cdot 2 \cdot 2 + 3 \cdot 2 + 1 \cdot 2 \cdot (n-4) = 3n+5$ if the added edges are adjacent, and $\xi^c(G) = (n-1) \cdot 1 + 2 \cdot 2 \cdot 4 + 1 \cdot 2 \cdot (n-5) = 3n+5$ if the added edges are not adjacent. This proves Corollary 3.

For a connected graph $G$, let $D_u = D_G(u) = \sum_{v \in V(G)} d(u,v|G)$. Then $D'(G) = \sum_{u \in V(G)} d_u D_u$ is the degree distance of $G$ [16,17], which is also a part of the Schultz molecular topological index [18-20]. We give a relation between the eccentric connectivity index and the degree distance.

**Proposition 4.** *Let $G$ be a connected graph with $n \geq 2$ vertices. Then*

$$\xi^c(G) \geq \frac{1}{n-1} D'(G)$$

*with equality if and only if $G = K_n$.*

**Proof.** Obviously, $e_u \geq \dfrac{D_u}{n-1}$ with equality if and only if $d(u,v|G)$ is a constant for all $v \in V(G)$ with $v \neq u$. It follows that

$$\xi^c(G) \geq \sum_{u \in V(G)} d_u \frac{D_u}{n-1} = \frac{1}{n-1} D'(G)$$

with equality if and only if $d(u,v|G)$ is a constant for all $u,v \in V(G)$ with $u \neq v$, i.e., $G = K_n$. □

Recall that the Wiener index of a connected graph $G$ is defined as $W(G) = \frac{1}{2}\sum_{u \in V(G)} D_u$ [21]. Let $G$ be a connected graph with $n \geq 2$ vertices and minimum degree $\delta$. By Proposition 4, it is easily seen that $\xi^c(G) \geq \frac{2\delta}{n-1}W(G)$ with equality if and only if $G = K_n$.

The first Zagreb index of a graph $G$ is defined as $M_1(G) = \sum_{u \in V(G)} d_u^2$ [22-25]. Let $K_n - ke$ be the graph formed by deleting $k$, where $k = 1,\ldots,\left\lfloor \frac{n}{2} \right\rfloor$, independent edges from the complete graph $K_n$. Let $K_n - 0e = K_n$. Obviously, $K_n - ke$ is actually a complete $(n-k)$-partite graph with exactly $n - 2k$ partite sets of cardinality one and $k$ partite sets of cardinality two, for $k = 0,1,\ldots,\left\lfloor \frac{n}{2} \right\rfloor$.

**Proposition 5.** *Let $G$ be a connected graph with $n \geq 3$ vertices and $m$ edges. Then*

$$\xi^c(G) \leq 2nm - M_1(G)$$

*with equality if and only if $G = K_n - ke$, for $k = 0,1,\ldots,\left\lfloor \frac{n}{2} \right\rfloor$, or $G = P_4$.*

**Proof.** Let $d(u;i)$ be the number of vertices that are of distance $i$ from a vertex $u$ in $G$, $i = 1,2,\ldots,e_u$. For $u \in V(G)$, it is easily seen that

$$n - 1 = d_u + \sum_{i=2}^{e_u} d(u;i) \geq d_u + \sum_{i=2}^{e_u} 1 = d_u + e_u - 1,$$

and thus $e_u \leq n - 1 - (d_u - 1) = n - d_u$ with equality if and only if $e_u = 1$ (and then $d_u = n-1$) or $e_u \geq 2$ with $d(u;2) = \cdots = d(u;e_u) = 1$. Then

$$\xi^c(G) = \sum_{u \in V(G)} d_u e_u \leq \sum_{u \in V(G)} d_u(n - d_u) = n \sum_{u \in V(G)} d_u - \sum_{u \in V(G)} d_u^2 = 2nm - M_1(G).$$

Suppose that equality holds in the above inequality. Then $e_u = n - d_u$, and thus $e_u = 1$ or $e_u \geq 2$ with $d(u;2) = \cdots = d(u;e_u) = 1$ for all $u \in V(G)$.

Suppose first that $e_u = 1$ for some $u \in V(G)$. Then $d_u = n - 1$ and $e_v = 1$ or 2 for all $v \neq u$. If $e_v = 1$ for all $v \neq u$, then $G = K_n$. Suppose that $e_v = 2$ for some $v \in V(G)$. Then there exists a vertex $w \in V(G)$ such that $d(v, w | G) = 2$. Since $d(v;2) = d(w;2) = 1$, the vertex $w$ is unique for fixed $v$. Thus, $d_v = d_w = n - 2$, implying that $G = K_n - ke$, for $k = 1, \ldots, \left\lfloor \frac{n-1}{2} \right\rfloor$.

Now suppose that $e_u \geq 2$ for all $u \in V(G)$. Then $d(u;2) = \cdots = d(u;e_u) = 1$. If $e_u = 2$ for all $u \in V(G)$, then $d_u = n - 2$ for all $u \in V(G)$, and since the sum of degrees $n(n-2)$ is always even, $n$ is even and thus $G = K_n - \frac{n}{2}e$. If $e_u \geq 3$ for some $u \in V(G)$, then $D(G) = 3$ (otherwise, for a center $z$ of a diametrical path, $d(z;2) \geq 2$, a contradiction), and then it is easily seen that $G = P_4$.

Conversely, it is easily checked that the upper bound for $\xi^c(G)$ is attained for $G = K_n - ke$, for $k = 0, 1, \ldots, \left\lfloor \frac{n}{2} \right\rfloor$, or $G = P_4$. □

Let $G$ be a connected graph with $n$ vertices, minimum degree $\delta$ and maximum degree $\Delta$. Then [26] $\sum_{u \in V(G)} e_u \leq \left( \frac{9n}{4\delta + 4} + \frac{15}{4} \right) n$, and thus $\xi^c(G) \leq \left( \frac{9n}{4\delta + 4} + \frac{15}{4} \right) \Delta n$.

## 4. RESULTS FOR TREES

In this section, we study the eccentric connectivity index for trees in more detail. We determine the trees of fixed diameter with the minimum eccentric connectivity index, and then deduce the $n$-vertex trees with respectively the minimum, second-minimum, and third-minimum eccentric connectivity indices for $n \geq 6$. We also determine the trees of fixed number of pendent vertices with the maximum eccentric connectivity index, and then deduce the $n$-vertex trees with respectively the maximum, second-maximum, and third-maximum eccentric connectivity indices for $n \geq 6$.

**Lemma 1.** *Let $u$ be a vertex of a tree $Q$ with at least two vertices. For integer $a \geq 1$, let $G_1$ be the tree obtained by attaching a star $S_{a+1}$ at its center $v$ to $u$ of $Q$, $G_2$ the tree obtained by attaching $a+1$ pendent vertices to $u$ of $Q$ (see Fig. 1). Then $\xi^c(G_2) < \xi^c(G_1)$.*

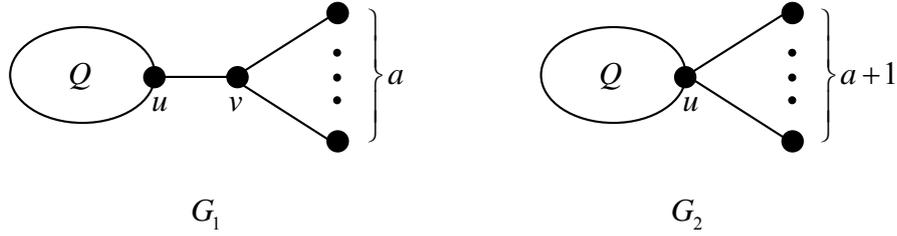

Fig. 1. The trees in Lemma 1.

**Proof.** Denote by $w$ a pendent neighbor of $v$ in $G_1$ and a pendent neighbor of $u$ in $G_2$ outside $Q$. It is easily seen that $e_{G_2}(x) \leq e_{G_1}(x)$ for $x \in V(Q)$, $e_{G_2}(w) = e_{G_1}(v)$ and $e_{G_1}(u) < e_{G_1}(w)$. Then

$$\xi^c(G_2) - \xi^c(G_1) = \sum_{x \in V(Q) \setminus \{u\}} d_Q(x) \left[ e_{G_2}(x) - e_{G_1}(x) \right]$$
$$+ \left[ d_Q(u) + a + 1 \right] e_{G_2}(u) - \left[ d_Q(u) + 1 \right] e_{G_1}(u)$$
$$+ 1 \cdot e_{G_2}(w) \cdot (a+1) - 1 \cdot e_{G_1}(w) \cdot a - (a+1) e_{G_1}(v)$$
$$\leq \left[ d_Q(u) + a + 1 \right] e_{G_1}(u) - \left[ d_Q(u) + 1 \right] e_{G_1}(u) - a e_{G_1}(w)$$
$$= a \left[ e_{G_1}(u) - e_{G_1}(w) \right] < 0,$$

and thus, $\xi^c(G_2) < \xi^c(G_1)$. The result follows. □

Let $\mathbf{T}(n,d)$ be the set of $n$-vertex trees with diameter $d$, where $2 \leq d \leq n-2$. Let $T_{(n,d)}$ be the set of $n$-vertex trees obtained from $P_{d+1} = v_0 v_1 \cdots v_d$ by attaching $n-d-1$ pendent vertices to $v_{\lfloor d/2 \rfloor}$ and/or $v_{\lceil d/2 \rceil}$, where $2 \leq d \leq n-2$. Note that in the set $T_{(n,d)}$, there is only one tree for even $d$, and $\left\lfloor \dfrac{n-d+1}{2} \right\rfloor$ trees for odd $d$. For a graph $G$ and a subset $E'$ of its edge set ($E^*$ of the edge set of its complement, respectively), $G-E'$ ($G+E^*$, respectively) denotes the graph formed from $G$ by deleting (adding, respectively) edges from $E'$ ($E^*$, respectively).

**Proposition 6.** *Let $G \in \mathbf{T}(n,d)$, where $2 \leq d \leq n-2$. Then*

$$\xi^c(G) \geq \left\lfloor \dfrac{3d^2+1}{2} \right\rfloor + (n-d-1)\left(1+2\left\lceil \dfrac{d}{2} \right\rceil\right)$$

*with equality if and only if $G \in \mathrm{T}_{(n,d)}$.*

**Proof.** The case $d=2$ is trivial.

Suppose that $d \geq 3$ and $G$ is a tree in $\mathbf{T}(n,d)$ with the minimum eccentric connectivity index. Let $P(G) = v_0 v_1 \cdots v_d$ be a diametrical path of $G$. By Lemma 1, all vertices outside $P(G)$ are pendent vertices adjacent to vertices of $P(G)$. Suppose that there exists some vertex $v_k$ with $k \neq \left\lfloor \dfrac{d}{2} \right\rfloor, \left\lceil \dfrac{d}{2} \right\rceil$, such that $d_G(v_k) \geq 3$. Suppose without loss of generality that $\left\lceil \dfrac{d}{2} \right\rceil < k \leq d-1$, i.e., $e_G(v_k) = k$. Denote by $w_1,\ldots,w_t$ by the pendent neighbors of $v_k$ outside $P(G)$. For $G_1 = G - \{v_k w_1, \ldots, v_k w_t\} + \{v_{\lfloor d/2 \rfloor} w_1, \ldots, v_{\lfloor d/2 \rfloor} w_t\} \in \mathbf{T}(n,d)$, we have

$$\begin{aligned}
\xi^c(G_1) - \xi^c(G) &= \left[ d_{G_1}(v_k)e_{G_1}(v_k) + 1 \cdot e_{G_1}(w_1) \cdot t + d_{G_1}(v_{\lfloor d/2 \rfloor})e_{G_1}(v_{\lfloor d/2 \rfloor}) \right] \\
&\quad - \left[ d_G(v_k)e_G(v_k) + 1 \cdot e_G(w_1) \cdot t + d_G(v_{\lfloor d/2 \rfloor})e_G(v_{\lfloor d/2 \rfloor}) \right] \\
&= \left[ 2k + t\left(1+\left\lceil \dfrac{d}{2} \right\rceil\right) + \left(d_G(v_{\lfloor d/2 \rfloor})+t\right)\left\lceil \dfrac{d}{2} \right\rceil \right] \\
&\quad - \left[ (t+2)k + t(1+k) + d_G(v_{\lfloor d/2 \rfloor})\left\lceil \dfrac{d}{2} \right\rceil \right] \\
&= 2t\left(\left\lceil \dfrac{d}{2} \right\rceil - k\right) < 0,
\end{aligned}$$

and thus, $\xi^c(G_1) < \xi^c(G)$, a contradiction. It follows that $d_G(v_i) = 2$ for all $1 \leq i \leq d-1$ with $i \neq \left\lfloor \frac{d}{2} \right\rfloor, \left\lceil \frac{d}{2} \right\rceil$. Thus $G \in T_{(n,d)}$.

Conversely, it is easily seen that $\xi^c(G) = \left\lfloor \frac{3d^2+1}{2} \right\rfloor + (n-d-1)\left(1 + 2\left\lceil \frac{d}{2} \right\rceil\right)$ for $G \in T_{(n,d)}$. □

For $d \geq 3$, it is easily checked that for $T' \in T_{(n,d)}$ and $T'' \in T_{(n,d-1)}$,

$$\xi^c(T') - \xi^c(T'') = \begin{cases} 2d-3 & \text{if } d \text{ is even} \\ 2n-3 & \text{if } d \text{ is odd} \end{cases}$$

and thus, $\xi^c(T') > \xi^c(T'')$. Note that $T_{(n,2)} = \{S_n\}$, $T_{(n,3)}$ contains exactly all the $\left\lfloor \frac{n-2}{2} \right\rfloor$ double-stars (formed by adding an edge connecting the centers of two stars, each with at least two vertices), and $T_{(n,4)}$ contains a unique tree formed by attaching $n-5$ pendent vertices to the center of the path with five vertices. Then, we have

**Proposition 7.** *Among the n-vertex trees with $n \geq 6$, $S_n$, n-vertex double-stars, and the tree formed by attaching $n-5$ pendent vertices to the center of the path with five vertices are respectively the unique trees with the minimum, second-minimum, and third-minimum eccentric connectivity indices, which are equal to $3(n-1)$, $5n-6$, and $5n-1$, respectively.*

**Lemma 2.** *Let $G_1$ and $G_2$ be the trees shown in Fig. 2, where vertices $x$ and $y$ are connected by a path of length at least one (if the length is more than one, then vertices in this path except $x$ and $y$ are of degree two), $y$ has a unique neighbor in $N$, and $d_{G_1}(x) \geq 2$. In $G_1$, $y$ has at least one neighbor in $M$, and all of such neighbors are switched to be neighbors of $x$ in $G_2$. If $\max\{d(y,u\mid G_1) : u \in V(M)\} \leq \max\{d(y,u\mid G_1) : u \in V(N)\}$ and $N$ is not a single vertex, then $\xi^c(G_1) < \xi^c(G_2)$.*

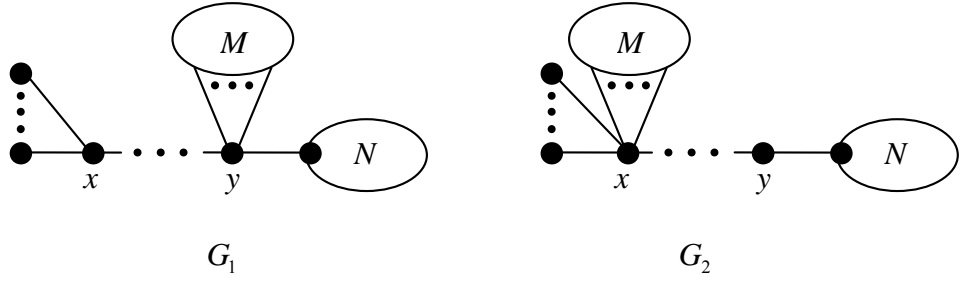

Fig. 2. The trees in Lemma 2.

**Proof.** Since $\max\{d(y,u|G_1): u \in V(M)\} \leq \max\{d(y,u|G_1): u \in V(N)\}$, we have $e_{G_2}(u) \geq e_{G_1}(u)$ for any $u \in V(G_1)$. Since $N$ is a not single vertex, we have $e_{G_1}(x) > e_{G_1}(y)$. Let $s$ be the number of neighbors of $y$ of $G_1$ in $M$. Then $d_{G_1}(y) = s+2$, $d_{G_2}(x) = d_{G_1}(x) + s$ and $s \geq 1$. We have

$$\begin{aligned}
\xi^c(G_2) - \xi^c(G_1) &= \sum_{u \in V(G_1)\setminus\{x,y\}} d_{G_1}(u)\left[e_{G_2}(u) - e_{G_1}(u)\right] \\
&\quad + d_{G_2}(x)e_{G_2}(x) - d_{G_1}(x)e_{G_1}(x) + d_{G_2}(y)e_{G_2}(y) - d_{G_1}(y)e_{G_1}(y) \\
&\geq \left[d_{G_1}(x) + s\right]e_{G_1}(x) - d_{G_1}(x)e_{G_1}(x) + 2e_{G_1}(y) - (s+2)e_{G_1}(y) \\
&= s\left[e_{G_1}(x) - e_{G_1}(y)\right] > 0,
\end{aligned}$$

which implies that $\xi^c(G_2) > \xi^c(G_1)$. □

Let $\mathscr{T}(n,p)$ be the set of $n$-vertex trees with $p$ pendent vertices, where $2 \leq p \leq n-2$. Let $T_a^{n,p}$ be the tree obtained by attaching $a$ and $p-a$ pendent vertices respectively to the two end vertices of the path $P_{n-p}$ for $1 \leq a \leq \lfloor \frac{p}{2} \rfloor$, and let $\mathrm{T}^{(n,p)} = \left\{T_a^{n,p} : 1 \leq a \leq \lfloor \frac{p}{2} \rfloor\right\}$. A path $u_1 u_2 \ldots u_r$ in a tree $T$ is said to be a pendent path rooted at $u_1$ if $d_T(u_1) \geq 3$, $d_T(u_i) = 2$ for $i = 2,\ldots,r-1$, and $d_T(u_r) = 1$. A pendent edge is a pendent path of length one.

**Proposition 8.** *Let $G \in \mathscr{T}(n,p)$, where $2 \leq p \leq n-2$. Then*

$$\xi^c(G) \leq \left\lfloor \frac{3(n-p+1)^2+1}{2} \right\rfloor + (p-2)(2n-2p+1)$$

with equality if and only if $G \in \mathrm{T}^{(n,p)}$.

**Proof.** The case $p=2$ is trivial. Suppose that $3 \leq p \leq n-2$ and $G$ is a tree in $\mathscr{T}(n,p)$ with the maximum eccentric connectivity index. There are some pendent paths (of length one or more) in $G$ as $p \geq 3$.

**Case 1.** There is some pendent path with length at least two in $G$.

Let $Q$ be a pendent path with length at least two in $G$. Let $u_1$ be the neighbor of the pendent vertex in $Q$ and $u_2$ the root in $Q$. Let $G'$ be the graph obtained from $G$ by deleting the vertices in $Q$. If $G'$ is not an empty graph, then making use of Lemma 2 to $G_1 = G$ by setting $x = u_1$ and $y = u_2$, we may get a tree in $\mathscr{T}(n,p)$ with larger eccentric connectivity index, a contradiction. Thus, $G'$ is an empty graph, i.e., $G = T_1^{n,p} \in \mathrm{T}^{(n,p)}$.

**Case 2.** Every pendent path in $G$ is a pendent edge.

In this case $p \geq 4$. Let $u$ be the neighbor of an end vertex and the neighbor of the other end vertex of a diameter-achieving path of $G$. Then both $u$ and $v$ have at least two pendent neighbors. Suppose that there exists a vertex on the path joining $u$ and $v$ with degree more than two. Let $w$ be such a vertex such that $d(u,w|G)$ is as small as possible. Making use of Lemma 2 to $G_1 = G$ by setting $x = u$ and $y = w$, we may get a tree in $\mathscr{T}(n,p)$ with larger eccentric connectivity index, a contradiction. Thus, either $d(u,v|G) = 1$, or the vertices in the path joining $u$ and $v$ are all of degree two, except $u$ and $v$. In either case, $G \in \mathrm{T}^{(n,p)}$.

Combining Cases 1 and 2, we have $G \in \mathrm{T}^{(n,p)}$. Conversely, it is easily seen that

$$\xi^c(G) = \left\lfloor \frac{3(n-p+1)^2+1}{2} \right\rfloor + (p-2)(2n-2p+1) \text{ for } G \in \mathrm{T}^{(n,p)}. \quad \square$$

Let $T_{n,i}$ be the tree obtained from $P_{n-1} = v_0 v_1 \cdots v_{n-2}$ by attaching a pendent vertex $v_{n-1}$ to $v_i$, where $1 \leq i \leq \left\lfloor \frac{n-2}{2} \right\rfloor$. Obviously, $\mathrm{T}^{(n,3)} = \{T_1^{n,3}\} = \{T_{n,1}\}$.

**Proposition 9.** *Among the n-vertex trees, where $n \geq 6$, $P_n$, $T_{n,1}$ and $T_{n,2}$ are respectively the unique trees with the maximum, second-maximum, and third-maximum eccentric connectivity indices, which are equal to $\left\lfloor \frac{3(n-1)^2+1}{2} \right\rfloor$, $\left\lfloor \frac{3(n-1)^2}{2} \right\rfloor - n$, and $\left\lfloor \frac{3(n-1)^2}{2} \right\rfloor - n - 2$, respectively.*

**Proof.** For $2 \leq p \leq n-3$, $T_1^{n,p} \in T^{(n,p)}$ and $T_1^{n,p+1} \in T^{(n,p+1)}$, by Lemma 1, we have $\xi^c(T_1^{n,p}) > \xi^c(T_1^{n,p+1})$. Note that $\mathscr{T}(n,2) = T^{(n,2)} = \{P_n\}$ and $T^{(n,3)} = \{T_{n,1}\}$. Then, by Proposition 8, $P_n$ and $T_{n,1}$ are respectively the unique $n$-vertex trees with the maximum and second-maximum eccentric connectivity indices.

Now suppose that $G$ is an $n$-vertex tree different from $P_n$ and $T_{n,1}$. Let $p$ be the number of pendent vertices of $G$.

If $p \geq 4$, then by the arguments as above, we have $\xi^c(G) \leq \xi^c(T^{(n,4)})$ with equality if and only if $G \in T^{(n,4)}$.

Suppose that $p = 3$. Then $G$ is a tree obtained by identifying three pendent vertices of three paths. Denote by $a, b$ and $c$ respectively the lengths of the three paths. Assume that $a \geq b \geq c$. Suppose first that $c = 1$. Then $G = T_{n,i}$ with $i \geq 2$. For $2 \leq i \leq \left\lfloor \frac{n-4}{2} \right\rfloor$, it is easily seen that

$$\xi^c(T_{n,i+1}) - \xi^c(T_{n,i})$$
$$= \left[ e_{T_{n,i+1}}(v_{n-1}) - e_{T_{n,i}}(v_{n-1}) \right] + \left[ 2e_{T_{n,i+1}}(v_i) - 3e_{T_{n,i}}(v_i) \right] + \left[ 3e_{T_{n,i+1}}(v_{i+1}) - 2e_{T_{n,i}}(v_{i+1}) \right]$$
$$= \left[ (n-2-i) - (n-1-i) \right] - (n-2-i) + (n-3-i) = -2 < 0,$$

and thus, $\xi^c(T_{n,\lfloor(n-2)/2\rfloor}) < \cdots < \xi^c(T_{n,2})$. Now suppose that $c \geq 2$. Let $x_1$ be the common vertex of the three paths, and $x_2$ the neighbor of the pendent vertex of the path with length $b$. Making use of Lemma 2 to $G_1 = G$ by setting $y = x_1$ and $x = x_2$, we have $\xi^c(G) < \xi^c(T_{n,c}) \leq \xi^c(T_{n,2})$. Hence, $T_{n,2}$ is the unique tree with the second-maximum eccentric connectivity index in $\mathscr{T}(n,3)$.

We are left to compare eccentric connectivity indices of $T_{n,2}$ and graphs in $T^{(n,4)}$. Note that $T_1^{n,4} = T_{n,2} - \{v_0 v_1\} + \{v_0 v_2\} \in T^{(n,4)}$. By Lemma 1, we have $\xi^c(T_1^{n,4}) < \xi^c(T_{n,2})$. Thus, $T_{n,2}$ is the unique $n$-vertex tree with the third-maximum eccentric connectivity index. □

## 5. CONCLUSION

In this paper, we deal with the eccentric connectivity index of a connected graph. We present various lower and upper bounds for the eccentric connectivity index in terms of other graph invariants including the number of vertices, the number of edges, the degree distance and the first Zagreb index, determine the $n$-vertex trees of diameter $d$, $3 \le d \le n-2$, with the minimum eccentric connectivity index, and the $n$-vertex trees of $p$ pendent vertices, $3 \le p \le n-2$, with the maximum eccentric connectivity index. In addition, we determine the $n$-vertex trees with the minimum, second-minimum and third-minimum eccentric connectivity indices as well as with the maximum, second-maximum and third-maximum eccentric connectivity indices for $n \ge 6$.

A necessary condition for a molecular descriptor to be an acceptable measure of branching is that within the set of all $n$-vertex trees, its values should be extremal for $P_n$ and $S_n$ [*e.g.*, 27]. Our result (in Propositions 7 and 9) shows that the eccentric connectivity index satisfies this requirement.

Some topics for further research may be: To determine the $n$-vertex tree(s) of given number of pendent vertices with the minimum eccentric connectivity index, the $n$-vertex tree(s) of given maximum degree with the maximum eccentric connectivity index, and the $n$-vertex connected graph(s) of $m$ edges, $n-1 \le m < \binom{n}{2}$, with the maximum eccentric connectivity index.

**Acknowledgement.** This work was supported by the Guangdong Provincial Natural Science Foundation of China (Grant No. 8151063101000026). The authors thank the referees for valuable comments and suggestions.